\documentclass[reqno,11pt]{article}

\usepackage{amsmath,latexsym,amssymb}

\input xy
\xyoption{all}

\def\sqr#1#2{{\vcenter{\hrule height.#2pt
        \hbox{\vrule width.#2pt height#1pt \kern#1pt
                \vrule width.#2pt}
        \hrule height.#2pt}}}
\def\square{\mathchoice\sqr64\sqr64\sqr{4}3\sqr{3}3}
\def\QED{\hfill$\square$\break}
\def\demo{\noindent{\bf Proof: }}

\def\geq{\geqslant}
\def\leq{\leqslant}

\def\ra{\rightarrow}
\def\lra{\longrightarrow}

\def\ext{\hbox{\rm Ext}}
\def\depth{\hbox{\rm depth}\,}
\def\proj{\hbox{\rm Proj}}
\def\spec{\hbox{\rm Spec}}

\def\shom{{\scriptstyle{\rm Hom}}}
\def\hom{\hbox{\rm Hom}}
\def\coker{\hbox{\rm coker}}
\def\homgr{\hbox{\rm Homgr}}

\def\adots{\mathinner{\mkern2mu\raise1pt\hbox{.}
\mkern3mu\raise4pt\hbox{.}\mkern1mu\raise7pt\hbox{.}}}

\def\om{{\omega}}
\def\im{{\frak m}}
\def\ib{{\frak b}}

\def\cu{{\cal C}}

\def\D{\Delta}

\newtheorem{Theorem}{\sc Theorem}[section]
\newtheorem{Lemma}[Theorem]{\sc Lemma}
\newtheorem{Corollary}[Theorem]{\sc Corollary}
\newtheorem{Proposition}[Theorem]{\sc Proposition}

\newtheorem{Remark}[Theorem]{\sc Remark}

\begin{document}

\baselineskip=13pt

\pagestyle{empty}

\ \vspace{1.7in}

\noindent {\LARGE\bf Liaison of varieties of small dimension and
  deficiency modules}

\vspace{.25in}

\noindent MARC \ CHARDIN, \ Institut de Math\'ematiques, CNRS \& Universit\'e
Paris VI, 4, place Jussieu, F-75252 Paris CEDEX 05,
 France.\hfill\break
\noindent {\it E-mail}: {\tt chardin@math.jussieu.fr}

\vspace{2.4cm}

\section{Introduction \hfill\break}

Liaison relates the cohomology of the ideal sheaf of a scheme to the
cohomology of the canonical module of its link. We here refer to
Gorenstein liaison in a projective space over a field: each ideal is
the residual of the other in one Gorenstein homogeneous
ideal of a polynomial ring. Assuming that the  
linked schemes (or equivalently one of them) are Cohen-Macaulay,
Serre duality expresses the cohomology of the canonical module in
terms of the cohomology of the ideal sheaf. Therefore, in the case of
Cohen-Macaulay linked schemes, the cohomology of ideal sheaves can
be computed one from another: up to shifts in ordinary and homological
degrees, they are exchanged and dualized. In terms of free
resolutions this means that, up to a degree shift, they may be obtained
one from another by dualizing the corresponding complexes (for
instance, the generators of one cohomology module corresponds to the
last syzygies of another cohomology module of the link). 

If the linked schemes are not Cohen-Macaulay, this property
fails. Nevertheless, experience on a computer shows that these
modules are closely related. We here investigate this relation
for surfaces and three-dimensionnal schemes. To describe our results,
note that the graded duals of the local cohomology modules are
Ext modules (into the polynomial ring or the Gorenstein quotient that
provides the linkage), let us set $\hbox{---}^{*}$ for the graded dual and
$D_{i}(M):=H^{i}_{\im}(M)^{*}$, the $i$-th deficiency module when
$i\not= \dim M$, while $\om_{M}=D_{\dim M}(M)$.  

In the case of surfaces, we have to understand the behaviour under
liaison of the Hartshorne-Rao module $D_{1}$ and the module
$D_{2}$. Note that $D_{2}$ has a finite part
$H^{0}_{\im}(D_{2})\simeq D_{0}(D_{2})^{*}$ and a quotient 
$D_{2}/H^{0}_{\im}(D_{2})\simeq D_{1}(D_{1}(D_{2}))$ that is either
$0$ or a Cohen-Macaulay module of dimension one. For a module $M$ of
dimension $d\geq 3$, there is an isomorphism
$H^{2}_{\im}(\om_{M})\simeq H^{0}_{\im}(D_{d-1}(M))$, which gives
$H^{2}_{\im}(\om )\simeq H^{0}_{\im}(D_{2})$ and, together with the
liaison sequence, shows how these three modules behave under linkage.     

In the case of dimension three, there are seven Cohen-Macaulay modules,
obtained by iterating Ext's on the deficiency modules, that encodes
all the information on the deficiency modules. Five of
these modules are permuted by liaison, up to duality and shift in
degrees, but this is not the case for the two others that sits in an
exact sequence with the corresponding two of the link. The key map in
this four terms sequence of finite length modules generalizes the
composition of the linkage map with Serre duality map and seems an
interesting map to investigate further.  

We begin with a short review on duality, and present some extensions of
Serre duality that we couldn't find in the literature. For instance if
$X$ is an equidimensional scheme such that 
$$
(*)\quad \depth {\cal O}_{X,x}\geq \dim {\cal O}_{X,x}-1,\quad \forall
x\in X,
$$
then the isomorphisms given by Serre duality (in the Cohen-Macaulay
case) are replaced by a long exact sequence (Corollary \ref{gensduality}) in which
one module over three is zero in the Cohen-Macaulay case.  In the
context of liaison, it turns out that condition $(*)$ is equivalent to
the fact that the $S_{2}$-ification of a link (or equivalently any
link) is a Cohen-Macaulay scheme (Proposition  \ref{S2-liaison}).

\bigskip

\section{Preliminaries on duality in a projective space \hfill\break}

Let $k$ be a commutative ring, $n$ a positive integer, $A:={\rm Sym}(k^{n})$,
$\im :=A_{>0}$ and $\om_{A/k}:=A[-n]$. Let
$\hbox{---}^{*}$ denote the graded dual into $k$ :  
$$
\hbox{---}^{*}:=\homgr_{A}(\hbox{---},k)=\bigoplus_{\mu\in{\bf
    Z}}\hom_{k}((\hbox{---})_{\mu},k)
$$ 
and  $\hbox{---}^{\vee}$ the dual into $\om_{A/k}$ : 
 $\hbox{---}^{\vee}:=\hom_{A} (\hbox{---},\om_{A/k})$.  

The homology of the \v Cech complex $\cu^{\bullet}_{\im}(\hbox{---})$ will be
denoted by $H^{\bullet}_{\im}(\hbox{---})$; notice that this homology is the
local cohomology supported in $\im$ (for any commutative ring $k$), as
$\im$ is generated by a regular sequence. \medskip

The key lemma is the following,

\begin{Lemma}\label{loc-duality}
If $L_{\bullet}$ is a graded complex of finitely
generated free $A$-modules, then 
\smallskip

{\rm (i)} $H^{i}_{\im}(L_{\bullet})=0$ for $i\not= n$,\smallskip

{\rm (ii)} there is a functorial  homogeneous isomorphism of degree zero:
$$
H^{n}_{\im}(L_{\bullet})\simeq (L_{\bullet}^{\vee})^{*}.
$$
\end{Lemma}
\demo Claim (i) is standard and (ii) directly follows from the classical
description of $H^{n}_{\im}(A)$ (see e.g. \cite[\S 62]{Se} or \cite[Ch. III,
Theorem 5.1]{Ha}). \QED
For an $A$-module $M$, we set
$D_{i}(M):=\ext^{n-i}_{A}(M,\om_{A/k})$ and
$$
\Gamma M:=\ker (\cu^{1}_{\im}(M)\ra \cu^{2}_{\im}(M))=\bigoplus_{\mu\in
  {\bf Z}}H^{0}({\bf P}^{n-1}_{k},\widetilde{M}(\mu )),
$$
so that we have an exact sequence
$$
0\ra H^{0}_{\im}(M)\ra M\ra \Gamma M\ra H^{1}_{\im}(M)\ra 0.
$$

We will also use the notation $\om_{M}:=D_{\dim M}(M)$ when $M$ is
finitely generated and $k$ is a field.  

The above lemma gives, 

\begin{Proposition}\label{ext-duality}
Assume that $M$ has a graded free resolution
$L_{\bullet}\ra M\ra 0$ such that each $L_{\ell}$ is finitely generated,
then for all $i$ there are functorial maps 
$$
H^{i}_{\im}(M)\buildrel{\phi^{i}}\over{\lra}
H_{n-i}((L_{\bullet}^{\vee})^{*})
\buildrel{\psi_{i}}\over{\lra} D_{i}(M)^{*}, 
$$
where $\phi^{i}$ is an isomorphism. Furthermore $\psi_{i}$ is an
isomorphism if $k$ is a field.
\end{Proposition}
\demo The two spectral sequences arising from the double
complex $\cu^{\bullet}_{\im}L_{\bullet}$ degenerates at the second
step (due to Lemma \ref{loc-duality} (i) for one, and the exactness of
localisation for the 
other) and together with Lemma \ref{loc-duality} (ii) provides the first
isomorphism. If $k$ is a field, $\hbox{---}^{*}$ is an exact functor,
so that the natural map
$H_{n-i}((L_{\bullet}^{\vee})^{*})\buildrel{\psi_{i}}\over{\lra} 
H^{n-i}(L_{\bullet}^{\vee})^{*}= D_{i}(M)^{*}$ is an isomorphism. \QED

\medskip
\begin{Proposition}\label{sseq-duality}
Assume that $M$ possses a
  graded   free resolution $L_{\bullet}\ra M\ra 0$ such that each
  $L_{\ell}$ is   finitely generated, then there is a spectral
  sequence 
$$
H^{i}_{\im}(D_{j}(M))\Rightarrow H_{i-j}(L_{\bullet}^{*}).
$$
In particular, if $k$ is a field, then there is a spectral sequence
$$
H^{i}_{\im}(D_{j}(M))\Rightarrow \left\{
\begin{matrix}M^{*}\quad &\hbox{if}\ i-j=0\\
  0\quad\quad &\hbox{else.}\\
\end{matrix}
\right.
$$
\end{Proposition}

\demo The two spectral sequences arising from the double
complex $\cu^{\bullet}_{\im}L_{\bullet}^{\vee}$ have as second terms
${'E}_{2}^{ij}=H^{i}_{\im}(\ext^{j}_{A}(M,\om_{A/k}))=H^{i}_{\im}(D_{n-j}(M))$,
${''E}_{2}^{ij}=0$ for $i\not= n$ and ${''E}_{2}^{nj}=
H^{j}(H^{n}_{\im}(L_{\bullet}^{\vee}))\simeq 
H^{j}((L_{\bullet}^{\vee\vee})^{*})\simeq
H^{j}(L_{\bullet}^{*})$. This gives the first result. If $k$ is a
field, $H^{j}(L_{\bullet}^{*})=H_{j}(L_{\bullet})^{*}$ which
gives the second claim.\QED

\begin{Corollary}\label{duality1}
If $k$ is a field, and $M$ is a finitely generated
graded $A$-module of dimension $d$, then there are functorial maps 
$$
d^{i}_{j}:={'d}_{2}^{i,n-j}:H^{i}_{\im}(D_{j}(M))\ra H^{i+2}_{\im}(D_{j+1}(M))
$$
such that\smallskip

{\rm (i)} $d^{0}_{d-1}:H^{0}_{\im}(D_{d-1}(M))\lra H^{2}_{\im}(\om_{M})$ is
an isomorphism if $d\geq 3$,

\smallskip
{\rm (ii)} if $d=3$ and $M$ is equidimensionnal and satisfies $S_{1}$, there
is an exact sequence, 
$$
0\ra H^{1}_{\im}(D_{2}(M))\buildrel{d^{1}_{2}}\over{\lra}
H^{3}_{\im}(\om_{M})\lra \Gamma M^{*}\lra 0,
$$

\smallskip
{\rm (iii)} if $d\geq 4$ there is a exact sequence,
$$
\xymatrix@C=25pt{
0\ar[r]& H^{1}_{\im}(D_{d-1}(M))\ar[r]^(.55){d^{1}_{d-1}}&
H^{3}_{\im}(\om_{M})\ar[r]^(.4){e}&
H^{0}_{\im}(D_{d-2}(M))&\ &\ \\}
$$
$$
\xymatrix@C=25pt{\ &\ 
\ar[r]^(.3){d^{0}_{d-2}}&
H^{2}_{\im}(D_{d-1}(M))\ar[r]^(.55){d^{2}_{d-1}}&H^{4}_{\im}(\om_{M})\ar[r]& C\ar[r]&0\\}
$$
where $e$ is the composed map
$$
e: 
\xymatrix@C=36pt{H^{3}_{\im}(\om_{M})\ar[r]^(.45){can}&\coker (d^{1}_{d-1})
\ar[r]^(.55){({'d}_{3}^{3,n-d})^{-1}}&\ker (d^{0}_{d-2})
\ar[r]^(.4){can}&H^{2}_{\im}(D_{d-2}(M))\\}
$$
such that  

{\rm (a)} if $d=4$ and $M$ is equidimensionnal and satisfies $S_{1}$,  $C$
sits in an exact sequence 
$$
0\ra H^{1}_{\im}(D_{2}(M))\buildrel{{'d}_{3}^{1,n-2}}\over{\lra}
C {\lra} \Gamma M^{*}\ra 0, 
$$

{\rm (b)} if $d\geq 5$, $C$ sits in an exact sequence
$$
0\ra \ker (d^{1}_{d-2})\ra C\ra \ker \left[ \ker (d^{0}_{d-3})\ra 
 \ker (d^{3}_{d-1})/\hbox{\rm im}(d^{1}_{d-2})\right]\ra 0.
$$
\end{Corollary}
\demo We may assume that $M$ does not have associated primes of
dimension $\leq 1$. Then these statements are direct consequences of
the spectral sequence of Proposition \ref{sseq-duality} applied to $\Gamma M$. Note
that if $M$ is equidimensional and satisfies $S_{1}$, then $\dim
D_{i}(M)<i$ for $i<d$ so that $H^{i}_{\im}(D_{i}(M))=0$ for $i\not=
d$.\QED

\begin{Remark}\label{duality1b} {\rm
One may dualize all the above maps into $k$. By
Proposition \ref{ext-duality} it gives rise to natural maps
$\theta_{ij}:D_{i+2}(D_{j+1}(M))\ra D_{i}(D_{j}(M))$ satisfying the
``dual statments'' of (i), (ii) and (iii). Note also that the map $e$
in (iii) gives a map $\epsilon :D_{0}(D_{d-2}(M))\lra
H^{0}_{\im}(D_{3}(\om_{M}))\subseteq D_{3}(\om_{M})$. Also in (iii)(a)
the graded dual $K$ of $C$ sits in an exact sequence
$$
0\ra \Gamma M\buildrel{\tau}\over{\lra} K\buildrel{\eta}\over{\lra}
D_{1}(D_{2}(M))\ra 0, 
$$
where $\tau$ is the graded dual over $k$ of the transgression map
in the spectral sequence and $\eta$ is the graded dual over $k$ of
${'d}_{3}^{1,n-2}$ (both composed with isomorphisms given by
Proposition \ref{ext-duality}).}
\end{Remark}

\medskip 

\begin{Corollary}
If $k$ is a field and $M$ is a finitely generated
equidimensionnal graded $A$-module of dimension $d$ such that
$\widetilde{M}$ satisfies $S_{\ell}$ for some $\ell\geq 1$, then, for
$1< i\leq \ell$, there are functorial surjective maps    
$$
f_{i}:H^{d+1-i}_{\im}(\om_{M})\lra D_{i}(M)=H^{0}_{\im}(D_{i}(M))
$$
which are isomorphismes for $1<i<\ell$. And there is 
an injection $D_{1}(M)\lra H^{d}_{\im}(\om_{M})$ whose cokernel is 
$(M/H^{0}_{\im}(M))^{*}$ if $\ell\geq 2$, so that in this case 
$H^{d}_{\im}(\om_{M})\simeq \Gamma M^{*}$.
\end{Corollary}
\demo The condition on $\widetilde{M}$ implies that
$H^{i}_{\im}(D_{j}(M))=0$ if $i>\max\{ 0,j-\ell \}$ and $j\not=
d$. These vanishings together with the spectral sequence of
Proposition \ref{sseq-duality} gives the result.\QED

The next result gives a slight generalization of Serre duality to
schemes (or sheaves) that are close to be Cohen-Macaulay:

\begin{Corollary}\label{gensduality}
 Let $k$ be a field and $M$ is a finitely generated
graded $A$-module which is unmixed of dimension $d\geq 3$. Set
$F_{M}:=D_{d-1}(M)$, ${\cal  M}:=\widetilde{M}$ and assume that
$$
\depth {\cal M}_{x}\geq \dim {\cal M}_{x}-1,
\quad \forall x\in \hbox{Supp}({\cal M}).
$$

Then the exists a functorial isomorphism $H^{0}_{\im}(F_{M})\simeq
H^{2}_{\im}(\om_{M})$ and  a long exact sequence 

$$
\xymatrix{
0\ar[r]&H^{1}_{\im}(F_{M})\ar[r]&H^{3}_{\im}(\om_{M})\ar[r]
&D_{d-2}(M)\ar[dll]&\\
&H^{2}_{\im}(F_{M})\ar[r]&H^{4}_{\im}(\om_{M})\ar[r]&
D_{d-3}(M)\ar[dll]&\\
&&\vdots &\ar[dll]&\\
&H^{d-3}_{\im}(F_{M})\ar[r]&H^{d-1}_{\im}(\om_{M})\ar[r]&
D_{2}(M)\ar[dll]\\
&H^{d-2}_{\im}(F_{M})\ar[r]&H^{d}_{\im}(\om_{M})\ar[r]&\Gamma M^{*}\ar[r]&0.\\}
$$
\end{Corollary}
\demo 
This is immediate from the spectral sequence of
Proposition \ref{sseq-duality}.\QED

\begin{Remark}\label{} {\rm
The vanishing of certain collections of modules $H^{i}_{\im}(D_{j}(M))$
corresponds to frequently used properties that $M$ may have. For
instance, assume that $M$ is equidimensionnal of dimension $d>0$ and
consider all the possibly non zero modules $H^{i}_{\im}(D_{j}(M))$ for
$j\not= d$:
$$
\xymatrix{
     &&\D_{d}\ar@{..}[dl]&\cdots &\D_{2}\ar@{..}[dl]&\D_{1}\ar@{..}[dl]\\
L_{d}\ar@{..}[r]&\ H^{0}_{\im}(D_{d-1})&\cdots&H^{0}_{\im}(D_{1})&H^{0}_{\im}(D_{0})
&\\
\vdots &\vdots &\adots & H^{1}_{\im}(D_{1})&& \\
L_{2}\ar@{..}[r]&\ H^{d-2}_{\im}(D_{d-1})&\adots&&&\\
L_{1}\ar@{..}[r]&\ H^{d-1}_{\im}(D_{d-1})&&&&\\
&&&&&\\
     &C_{d}\ar@{..}[u]&\cdots &C_{2}\ar@{..}[uuu]&C_{1}\ar@{..}[uuuu]&\\}
$$
Then,

(i) The modules on lines $L_{1},\ldots ,L_{\ell}$ are zero if and only
if $M$ is Cohen-Macaulay in codimension $\ell$,\smallskip

(ii) The modules on columns $C_{1},\ldots ,C_{\ell}$ are zero if and only
if $M$ has depth at least $\ell$, 

(iii) The modules on diagonals $\D_{1},\ldots ,\D_{\ell}$ are zero if and only
if $M$ satisfies $S_{\ell}$.\smallskip

Note also that the spectral sequence shows that $\om_{M}$ is
Cohen-Macaulay if all non zero modules in this diagram are located on
diagonals $\D_{1}$ and $\D_{2}$ (even if $M$ is not
equidimensionnal). This generalizes a little the case of a
sequentially Cohen-Macaulay module (due to Schenzel), where all non
zero modules are located on diagonal $\D_{1}$. 
}

\end{Remark}

\section{The liaison sequence \hfill\break}

In this paragraph $A$ is a polynomial
ring over a field $k$. Recall that an homogeneous ideal $\ib$ of
$A$  is Gorenstein if $B:=A/\ib$ is Cohen-Macaulay and $\om_{B}$ is a
free $B$-module of rank 1 ({\it i.e.} $B$ is Gorenstein), so that
$\om_{B}\simeq B[a]$ for some $a\in {\bf Z}$ called the $a$-invariant
of $B$. The arithmetically Gorenstein subschemes of $\proj (A)$ are
the subschemes of the form $\proj (A/\ib )$ for an homogeneous
Gorenstein ideal $\ib$; if the scheme is not empty, such a $\ib$ is
unique.  
\medskip

\begin{Proposition}\label{liaison}
Let $I$ and $J=\ib:I$ be two ideals of $A$, linked
by a Gorenstein homogeneous ideal $\ib$ of $A$ (so that $I=\ib
:J$). Set $B:=A/\ib $, $R:=A/I$, $S:=A/J$, $d:=\dim B$, 
$F_{S}:=D_{d-1}(S)$ and $E_{R}:={\rm  End}(\om_{R})$. There are
natural exact sequences or isomorphisms:
$$
\leqno{\rm (i)}\quad\xymatrix@C=25pt{
0\ar[r]&\om_{R}\otimes \om_{B}^{-1}\ar[r]^(.4){\iota}
&\om_{B}\otimes_{B}\om_{B}^{-1}=B\ar[r]^(.7){s}&S\ar[r]&0\\
}
$$
where $\hbox{---}^{-1}:=\hom_{B}(\hbox{---},B)$ and $\iota$ and $s$
are the canonical maps; 
$$
\leqno{\rm (ii)}\quad\xymatrix@C=55pt{
0\ra \om_{S}\ar[r]^(.55){\shom_{B}(s,1_{\om_{B}})}&\om_{B}
\ar[r]^(.45){\shom_{B}(\iota,1_{\om_{B}})}&E_{R}\otimes_{B}\\}  
\om_{B}\buildrel{\delta^{r}}\over{\lra}F_{S}\ra 0
$$
where $\delta^{r}$ is the connecting map in the Ext sequence derived
from (i), in particular $F_{S}\simeq (E_{R}/R)\otimes_{B}\om_{B}$;
$$
\leqno{\rm (iii)}\quad \lambda_{i}:D_{i}(S)\buildrel{\sim}\over{\lra}
D_{i+1}(\om_{R})\otimes_{B}\om_{B} 
\quad \forall i\leq d-2
$$
given by the connecting map in the Ext sequence derived
from (i).
\end{Proposition}

\demo (i) is standard ({\it see} e.g. \cite{PS} or \cite{Mi}). Also (ii) and (iii)
directly follows by the Ext sequence derived from (i).  \QED

\begin{Remark} {\rm
(i) For any $B$-module $M$, there are functorial isomorphisms: 
$$
\tau_{i}:D_{i}(M)\buildrel{\sim}\over{\lra} \ext^{d-i}_{B}(M,\om_{B}).
$$

(ii) Choosing an isomomorphism $\phi: B[-a]\lra \om_{B}^{-1}$, (1)
gives an exact sequence,  
$$
\xymatrix@C=40pt{
0\ar[r]&\om_{R}[-a]\ar[r]^(.55){\iota\circ(1\otimes\phi )}
&B\ar[r]^{s}&S\ar[r]&0\\
}
$$
}
\end{Remark}

\begin{Proposition}\label{S2-liaison}
Let $X$ and $Y$ be two schemes linked by an
arithmetically Gorenstein projective subcheme of $\proj (A)$ of
dimension $d\geq 1$ and $a$-invariant $a$. Let  $Y_{2}:=\spec ({\cal
  E}nd (\om_{Y}))$ be the $S_{2}$-ification of $Y$. The following are
equivalent:\smallskip 

{\rm (i)} $\depth {\cal O}_{X,x}\geq \dim  {\cal O}_{X,x}-1$ for every $x\in
X$,\smallskip

{\rm (ii)} $\om_{Y}$ is Cohen-Macaulay,\smallskip

{\rm (iii)} $Y_{2}$ is Cohen-Macaulay,\smallskip

{\rm (iv)} $H^{i}(X,{\cal O}_{X}(\mu ))\simeq
H^{d-i-1}(Y,{\cal O}_{Y_{2}}(a-\mu ))^{*}$, for every $0<i<d-1$ and
every $\mu$, where $\hbox{---}^{*}$ denotes the dual into $k$. 
\end{Proposition}

\demo We may assume $d\geq 3$. The equivalence of (i) and (ii)
follows from Proposition  \ref{liaison}(iii). Now ${\cal
  O}_{Y_{2}}=\om_{\om_{Y}}$ and  $\om_{Y}=\om_{Y_{2}}$ which proves
that (ii) and (iii) are  equivalent.  

On the other hand, if $X$ satisfies (i), Corollary 
\ref{gensduality} and Proposition \ref{liaison} (ii) and (iii) implies (iv). 

If (iv) is satisfied, Proposition \ref{liaison}(iii) and the equality 
$\om_{Y}=\om_{Y_{2}}$ shows that $H^{i}(Y,{\cal O}_{Y_{2}}(\mu
))\simeq H^{d-i}(Y,\om_{Y_{2}}(-\mu ))^{*}$ for $i\not= d-1$. Applying
Corollary \ref{gensduality} with $M:=\Gamma_{*}{\cal O}_{Y_{2}}$ then implies that
$F_{M}$ has finite length, and therefore (iii) is satisfied. \QED 

\begin{Remark} {\rm
If one of the equivalent conditions of
Proposition \ref{S2-liaison} is satisfied, and $R$ and $S$ are the standard
graded unmixed algebras defining respectively $X$ and $Y$, then 
$H^{d-1}(X,{\cal O}_{X}(\mu ))\simeq H^{d}_{\im}(R)_{\mu}\simeq D_{d}(R)_{-\mu}$ and
$H^{0}_{\im}(D_{d}(R))\simeq H^{2}_{\im}(\om_{R})\simeq H^{1}_{\im}(S)[a]$.}
\end{Remark}

\section{Cohomology of linked surfaces and three-folds \hfill\break}

Let $k$ be a field, $A$ a polynomial ring over $k$ and
$\om_{A}:=A[-\dim A]$.\medskip

As in the first section, $\hbox{---}^{*}$ will denote the graded dual into
$k$, we set $D_{i}(M):=\ext^{\dim   A-i}_{A}(M,\om_{A})$ and
we introduce for simplicity the following abreviated notation : $D_{ij\cdots
  l}(M):=D_{i}(D_{j}(\cdots (D_{l}(M))))$.  

\subsection{The surface case \hfill\break}

Let $R$ be a homogeneous quotient of $A$
defining a surface 
$X\subset\proj (A)$. We will asume that $H^{0}_{\im}(R)=0$, as
replacing $R$ by $R/H^{0}_{\im}(R)$ leave both $X$ and
$H^{i}_{\im}(R)$ for $i>0$ unchanged. Then 
the three Cohen-Macaulay modules $D_{1}(R)$,  $D_{02}(R)$ and
$D_{12}(R)$ sits in the following exact sequences that shows how they
encode all the non-ACM character of $X$:
$$
0\lra R\lra \Gamma R \lra D_{1}(R)^{*}\lra 0
$$
$$
0\lra D_{02}(R)^{*}\lra D_{2}(R)\lra \Gamma D_{2}(R) \lra D_{12}(R)^{*}\lra 0
$$
where $\Gamma D_{2}(R)=\oplus_{\mu\in {\bf Z}}V_{\mu}$ and $\delta :=\dim
V_{\mu}$ is a constant which is zero if and only if $X$ is
Cohen-Macaulay. Also $D_{02}(R)$ is  a finite length module and $D_{12}(R)$
is either zero or a Cohen-Macaulay module of dimension one and degree
$\delta$.\medskip

The next result shows how these Cohen-Macaulay modules that encodes the
non-ACM character of $X$ behave under linkage.

\begin{Proposition}\label{surfaces}
Let $X=\proj (R)$ and $Y=\proj (S)$ be two
  surfaces linked by an arithmetically Gorenstein subcheme
  $\proj (B)\subseteq \proj (A)$ so that $\om_{B}\simeq B[a]$. Then
  there are natural degree zero isomorphisms,
$$
D_{1}(S)\simeq D_{02}(R)\otimes_{B}\om_{B},\ \ 
D_{12}(S)\simeq D_{112}(R)\otimes_{B}\om_{B}^{-1}.
$$
\end{Proposition}

Note that the roles of $R$ and $S$ may be reversed, so that
$D_{02}(S)\simeq D_{1}(R)\otimes_{B}\om_{B}^{-1}$ and
$D_{112}(S)=D_{12}(R)\otimes_{B}\om_{B}$.\medskip 

\demo By Proposition \ref{liaison} (iii), $D_{1}(S)\simeq
D_{2}(\om_{R})\otimes_{B}\om_{B}$ and $D_{2}(\om_{R})\simeq D_{02}(R)$
by  Corollary \ref{duality1} (i), this proves the first claim.

As for the second claim,  Corollary \ref{duality1} (ii) gives an exact sequence
$$
0\lra \Gamma S\lra D_{3}(\om_{S})\lra D_{12}(S)\lra 0
$$
where $D_{3}(\om_{S})\simeq E_{S}$. Now $\Gamma S/S\simeq
D_{1}(S)^{*}$ so that we have an exact sequence, 
$$
0\lra D_{1}(S)^{*}\lra E_{S}/S \lra D_{12}(S)\lra 0
$$
but $E_{S}/S\simeq D_{2}(R)\otimes_{B}\om_{B}^{-1}$ by Proposition \ref{liaison}
(ii) and we know that $D_{1}(S)^{*}\simeq
D_{02}(R)^{*}\otimes_{B}\om_{B}^{-1}$ so that the above sequence gives
$$
0\lra D_{02}(R)^{*}\otimes_{B}\om_{B}^{-1}\lra
D_{2}(R)\otimes_{B}\om_{B}^{-1}\lra D_{12}(S)\lra 0 
$$
but  $D_{02}(R)^{*}=H^{0}_{\im}(D_{2}(R))$ and
$D_{112}(R)=D_{2}(R)/H^{0}_{\im}(D_{2}(R))$.\QED

\subsection{The three-fold case \hfill\break}

Let $R$ be a homogeneous quotient of $A$
defining a three-fold $X\subset\proj (A)$. As in the surface case will
asume that $H^{0}_{\im}(R)=0$. Now the seven Cohen-Macaulay modules
$D_{1}(R)$, $D_{02}(R)$ and 
$D_{12}(R)$, $D_{03}(R)$, $D_{013}(R)$, $D_{113}(R)$ and $D_{23}(R)$
also encode all the non-ACM character of $X$.
\medskip

The following result essentially gives the behaviour of these modules
under liaison.

\begin{Proposition}\label{three-fold}
Let $X=\proj (R)$ and $Y=\proj (S)$ be two
  three-folds linked by an arithmetically Gorenstein subcheme
  $\proj (B)\subseteq \proj (A)$ so that $\om_{B}\simeq B[a]$. Then
  there are natural degree zero isomorphisms,
$$
D_{1}(S)\simeq D_{03}(R)\otimes_{B}\om_{B},
D_{12}(S)\simeq D_{113}(R)\otimes_{B}\om_{B}^{-1},
D_{23}(S)\simeq D_{223}(R)\otimes_{B}\om_{B}^{-1}.
$$
And there is an exact sequence,
$$
\xymatrix@C=18pt{
0\ar[r]&D_{013}(R)\ar[r]&D_{02}(S)\otimes_{B}\om_{B}\ar[r]^(.58){\psi}&D_{002}(R)\ar[r]
&D_{0013}(S)\otimes_{B}\om_{B}\ar[r]&0.\\ 
}
$$
where $\psi$ is the composition of the maps from Remark  \ref{duality1b} and
Proposition \ref{liaison} (iii):
$$
\xymatrix@C=50pt{
D_{02}(S)\otimes_{B}\om_{B}\ar[r]^(.55){D_{0}(\lambda_{2}^{-1})\otimes
  1_{\om_{B}}}&D_{03}(\om_{R})\ar[r]^{\epsilon^{*}}&D_{002}(R).\\}
$$
\end{Proposition}

Note that the roles of $R$ and $S$ may be reversed, so that we have a
complete list of relations between these fourteen modules. Also the
last exact sequence may be written
$$
\xymatrix@C=16pt{
0\ar[r]&D_{013}(R)\ar[r]&D_{02}(S)\otimes_{B}\om_{B}\ar[r]^(.58){\psi}&D_{02}(R)^{*}\ar[r]
&D_{013}(S)^{*}\otimes_{B}\om_{B}\ar[r]&0.\\ 
}
$$

\demo $D_{1}(S)\simeq D_{2}(\om_{R})\otimes_{B}\om_{B}\simeq
D_{03}(R)\otimes_{B}\om_{B}$ by  Proposition \ref{liaison} (iii) and Corollary 
\ref{duality1} (i). Also  Corollary \ref{duality1} (iii) provides exact sequences:
$$
\hbox{(1)}\quad
\xymatrix@C=15pt{
0\ar[r]& K\ar[r]^{can}& E_{R}\ar[r]^(.4){\theta_{23}}&D_{23}(R)
\ar[r]^{\theta_{02}}&D_{02}(R)\ar[r]^{e}&D_{3}(\om_{R})\ar[r]^{\theta_{13}}& 
D_{13}(R)\ar[r]& 0\\}
$$
and 
$$
\hbox{(2)}\quad \xymatrix{0\ar[r]& \Gamma R\ar[r]^{\tau}&
  K\ar[r]^(.4){\eta}& D_{12}(R)\ar[r]&0\\}. 
$$
From (2) we get isomorphisms $D_{i}(\tau ):D_{i}(K)\buildrel{\sim
}\over{\lra}D_{i}(R)$ for $i\geq 3$ and another exact sequence 
$$
\hbox{(3)}\quad\xymatrix{
0\ar[r]&D_{2}(K)\ar[r]^{D_{2}(\tau )}&D_{2}(R)\ar[r]^{\delta}&
D_{112}(R)\ar[r]^{D_{1}(\eta )}&D_{1}(K)\ar[r]& 0\\},
$$
from which it follows that $D_{1}(K)=0$ (note that this vanishing also
follows upon taking Ext's on sequence (1)).
As $D_{02}(R)$ is of dimension zero, the right end of the first
sequence then shows that
$D_{1}(\theta_{13}):D_{113}(R)\buildrel{\sim}\over{\lra}D_{13}(\om_{R})$,
but we have an isomorphism $\lambda_{2}^{-1}\otimes 1_{\om_{B}^{-1}}:D_{3}(\om_{R}) \buildrel{\sim}\over{\lra}D_{2}(S)\otimes_{B}\om_{B}^{-1}$ by
Proposition 
\ref{liaison} (iii) and this provides the second isomorphism. 
Now taking Ext's on (1) also provides two exact sequences 
$$
\hbox{(4)}\quad \xymatrix{
0\ar[r]&D_{013}(R)\ar[r]^{D_{0}(\theta_{13})}&D_{03}(\om_{R})\ar[r]^{\epsilon^{*}}&D_{002}(R)\ar[r]^(.6){can}&L\ar[r]& 0\\} 
$$
$$
\hbox{(5)}\ \xymatrix@C=14pt{
0\ar[r]&D_{3}(E_{R})\ar[r]& D_{3}(R)\ar[r]& D_{223}(R)\ar[r]& D_{2}(E_{R})\ar[r]&
D_{2}(K)\ar[r]^(.6){can}& L'\ar[r]& 0\\}.
$$
and an isomorphism $\mu :L\ra L'$. Now (5) together with (3) gives us a complex
$$
\hbox{(6)}\ \xymatrix@C=10pt{
0\ar[r]&D_{3}(E_{R})\ar[r]& D_{3}(R)\ar[r]& D_{223}(R)\ar[r]& D_{2}(E_{R})\ar[r]&
D_{2}(R)\ar[r]& D_{112}(R)\ar[r]& 0\\}
$$
whose only homology is the subquotient $L$ of $D_{2}(R)$. 

Taking Ext's on the exact sequence $0\ra R\ra E_{R}\ra
D_{3}(S)\otimes_{B}\om_{B}^{-1}\ra 0$ given by Propostion \ref{liaison} (ii), we get
an exact sequence
$$
\hbox{(7)}\quad 0\ra D_{3}(E_{R})\ra D_{3}(R)\ra
D_{23}(S)\otimes_{B}\om_{B}\ra D_{2}(E_{R})\quad\quad \quad \quad \quad  
$$
$$
\quad \quad \quad \quad \quad \quad \quad \quad \quad \quad \quad \quad 
\ra D_{2}(R)\ra D_{13}(S)\otimes_{B}\om_{B}\ra 0.
$$
Comparing (6) and (7) proves the third isomorphism and provides the
exact sequence
$$
\hbox{(8)}\quad 0\ra L\ra D_{13}(S)\otimes_{B}\om_{B}\ra D_{112}(R)\ra
0
$$
showing that $L\simeq  D_{0013}(S)\otimes_{B}\om_{B}$ because
$D_{112}(R)\simeq D_{1113}(S)\otimes_{B}\om_{B}$. This last
identification together with (4) provides the exact sequence of the
Proposition and concludes the proof.\QED

\smallskip


\end{document}